\documentclass[12pt]{amsart}
\usepackage[margin=1in]{geometry}
\usepackage{amsmath,amssymb,amsthm}
\usepackage{hyperref}
\usepackage{cleveref}

\newtheorem{theorem}{Theorem}[section]
\newtheorem{lemma}[theorem]{Lemma}

\newtheorem{proposition}[theorem]{Proposition}
\theoremstyle{definition}
\newtheorem{definition}[theorem]{Definition}

\theoremstyle{remark}
\newtheorem{remark}[theorem]{Remark}

\newcommand{\R}{\mathbb R}
\newcommand{\C}{\mathbb C}
\newcommand{\N}{\mathbb N}
\newcommand{\Hn}{\mathbb H^n}
\newcommand{\FS}{\mathcal F_S}
\newcommand{\WS}{W_{S,\lambda}}
\newcommand{\phik}{\varphi^{n-1}_{k,\lambda}}
\newcommand{\Lk}{L^{n-1}_k(0)}

\title[A Characterisation of the Strichartz Fourier Transform] {A Characterisation of the Strichartz Fourier Transform
via the Strichartz Weyl Transform}
\author{R. Lakshmi Lavanya}
\address{Department of Mathematics, Indian Institute of Science Education and Research (IISER) Tirupati, Yerpedu, Tirupati 517619, Andhra Pradesh, India.}
\email{rlakshmilavanya@labs.iisertirupati.ac.in}

\date{}

\begin{document}
\maketitle

\begin{abstract}
Thangavelu recently introduced the Strichartz Fourier transform, a scalar-valued Fourier transform on the Heisenberg group. Lakshmi Lavanya and Thangavelu had earlier characterised the Weyl transform, as essentially the unique continuous *-homomorphism of the twisted convolution algebra into the bounded operators on square-integrable functions, and used this to characterise the group Fourier transform on the Heisenberg group. Here we obtain an analogous pair of results for the Strichartz Fourier transform. Because this transform is scalar-valued, it does not factor through an operator-valued object of the kind used in that case. The role of the Weyl transform is instead played by the Strichartz Weyl transform, a recently introduced scalar-valued analogue. We characterise the Strichartz Weyl transform by four properties: a convolution-reproduction property at a single test function, a twisted translation-covariance property, a unitary-invariance property, and a single normalisation. We then deduce from this a characterisation of the Strichartz Fourier transform itself on each ray of the Heisenberg fan. At the boundary ray of the fan, the twisted structure degenerates. There, the analogous characterisation is proved by Euclidean methods.
\end{abstract}



\section{Introduction}

The Fourier transform on a locally compact abelian group $G$ can be characterised purely by its algebraic behaviour under convolution. Jaming \cite{Jaming} showed that a continuous linear map $T:L^1(G)\to C(\widehat G)$ satisfying $T(f*g)=T(f)T(g)$ for all $f,g\in L^1(G)$ must be of the form
\begin{align*}
T(f)(\xi)=\chi_E(\xi)\,\widehat f(\varphi(\xi))
\end{align*}
for a set $E\subset\widehat G$ and a suitable map $\varphi$. In particular, the Fourier transform is the essentially unique continuous multiplicative linear map into $C(\widehat G)$. Jaming's characterisation was itself preceded by Alesker, Artstein-Avidan, and Milman's characterisation of the Fourier transform on $\R^n$ \cite{AAM}. The two results are quite different in nature. Later, many characterisations of the Fourier transform were obtained in \cite{K1,K2,KS}.

On the Heisenberg group $\Hn$, the group Fourier transform $\widehat f(\lambda)$ is operator-valued, taking values in $\mathcal B(L^2(\R^n))$. Jaming's result does not apply directly. Lakshmi Lavanya and Thangavelu \cite{LLT} proved the natural analogue. A nonzero continuous linear map $T$ satisfying the adjoint property, the operator-valued convolution property
\begin{align*}
T(f*g)(\lambda)=T(f)(\lambda)\,T(g)(\lambda),
\end{align*}
and the translation-covariance property
\begin{align*}
T(R_{(z,t)}f)(\lambda)=T(f)(\lambda)\,\pi_\lambda(z,t)^*,
\end{align*}
must coincide with $\widehat f(\lambda)$ at every $\lambda$ at which it does not vanish identically. 
Vector-valued analogues of characterisations of the Fourier transform on various domains have been obtained in \cite{ABM1,ABM2,AM}.

Recently, Thangavelu \cite{Th23} introduced a different Fourier transform on $\Hn$, namely, the Strichartz Fourier transform. It is a scalar-valued transform built from the joint spectral theory of the sublaplacian $\mathcal L$ and the central vector field $\partial/\partial t$ of the Heisenberg group, rather than from the representation theory of $\Hn$. This transform shares several structural features with the Helgason Fourier transform on rank one Riemannian symmetric spaces. Its scalar-valuedness resolves technical difficulties present for the operator-valued transform. These arise in describing the image of the Schwartz space and in formulating sharp weighted inequalities. Since its introduction, $\FS$ has already served as the basis for further work. Recent papers have established uncertainty principles \cite{DDG} and weighted norm inequalities of Pitt type \cite{DGPT} for $\FS$.

More recently, a characterisation of $\FS$ along a different route has been obtained by Lavanya (\cite{Ksft}), via a differentiation/shift-property framework adapted from Krukowski's work on the Euclidean and locally compact abelian settings. There, $\FS$ is determined on each ray by a joint spectral property, formulated directly through the action of the sublaplacian $\mathcal L$ and the central vector field $T$, together with a twisted translation property and a normalisation at a single test function (and, alternatively, by the twisted translation property alone paired with a Dirac-delta-type approximate-identity hypothesis). The approach taken here is of a different character. $W_{S,\lambda}$ is characterised by its behaviour under the twisted convolution product itself: a convolution-reproduction property at a single test function, together with twisted translation covariance and $U(n)$-invariance. This mirrors the convolution-algebra strategy of \cite{LLT} rather than the differential one of \cite{Khn,Ksft}. The two approaches reach $\FS$ by essentially disjoint mechanisms, and neither set of hypotheses is a consequence of the other.

Because $\FS$ is scalar-valued, however, it cannot factor through a faithful representation of the twisted convolution algebra. There is no room, in $C(\C^n)$, for a matrix-unit relation of the operator-theoretic kind that drives the second stage of Lavanya and Thangavelu's argument \cite{LLT}. We mirror their broad two-stage strategy. First we characterise the $\C^n$-level object, then we lift it via translation covariance. But the technique carrying out each stage here is necessarily of a different kind, and the resulting object satisfies only a convolution property restricted to a single test function. This raises the question of what plays the role of the Weyl transform for $\FS$.

Central translation covariance identifies the relevant reduced object as in the classical case, up to the sign convention. With the covariance of Theorem~\ref{thm:FSchar}, which carries the factor $e^{+i\lambda\tau}$, $\FS(f)(a,z)$ depends on $f$ only through $f^{-\lambda}$, for $a=(\lambda,(2k+n)|\lambda|)$ on the $k$-th ray of the Heisenberg fan. Explicitly,
\begin{align*}
\FS(f)(a,z)=f^{-\lambda}*_{-\lambda}\phik(z)
\end{align*}
(notation as in Section~\ref{sec:prelim}). Here $\phik$ is the Laguerre function of type $n-1$ and degree $k$. The map $g\mapsto g*_{-\lambda}\phik$ is, up to the scalar $(2\pi)^n|\lambda|^{-n}$, the $k$-th spectral projection $P_k$ of the special Hermite operator $L_{-\lambda}$ on $L^2(\C^n)$. Its properties as a projection, in relation to twisted convolution, and in relation to the special Hermite expansion are well established in the literature on the Heisenberg group and its Gelfand pairs. This map plays precisely the structural role that the Weyl transform plays for the Fourier transform. It was identified and named the Strichartz Weyl transform in \cite{SWTpaper}. Beyond the elementary properties recalled in Section~\ref{sec:prelim} below, that paper develops several further properties of this transform. None of this further material is needed for, or used in, the characterisation proved here. We obtain a structural characterisation of $\WS$ itself, and then deduce one for $\FS$.

Since $\WS$ is a projection and not a representation, the matrix-unit technique that closes the argument is unavailable here. Hypothesis (i) below asks only that the single Laguerre function $\phik$ be reproduced correctly under twisted convolution. On its own it says nothing about how $S$ acts on any other function, until a convolution-kernel representation for $S$ has first been produced by other means. We show in Section~\ref{sec:charWS} that two further properties, a twisted translation-covariance property on $\C^n$ and a $U(n)$-invariance property, first manufacture such a kernel and force it to be radial. Only then does the single-function hypothesis identify the kernel as $\phik$ itself. This gives a complete characterisation of $\WS$ (Theorem~\ref{thm:char}). We then lift this result, in Section~\ref{sec:charFS}, via full translation covariance on $\Hn$, to a characterisation of $\FS$ itself (Theorem~\ref{thm:FSchar}). The $\C^n$-level twisted translation is the restriction of this full translation covariance.

The characterisation just described covers the rays $R_k$, $k\in\N$, of the Heisenberg fan $\Omega$. It says nothing about the boundary ray $R_\infty=\{(0,\tau):\tau\ge0\}$, where $\lambda=0$ and the twisted structure underlying Sections~\ref{sec:charWS}-\ref{sec:charFS} degenerates. There, $\WS$ and $\FS$ reduce, respectively, to the spectral density of the ordinary Laplacian on $\C^n\cong\R^{2n}$ and to the Hankel transform. Section~\ref{sec:boundary} records this degeneration explicitly. It proves the corresponding characterisation at $R_\infty$ by methods of Euclidean harmonic analysis rather than the representation-theoretic argument used elsewhere in the paper, so that $\FS$ is characterised at each point of $\Omega$.

The paper is organised as follows. Section~\ref{sec:prelim} recalls the necessary notations and discusses the properties used in the subsequent sections. Section~\ref{sec:charWS} proves the characterisation of $\WS$, along with an alternative proof via the ellipticity of the special Hermite operator. Section~\ref{sec:charFS} verifies that $\FS$ satisfies the corresponding hypotheses on $\Hn$ and proves the characterisation of $\FS$, reducing it to Section~\ref{sec:charWS}. Section~\ref{sec:boundary} treats the boundary ray $R_\infty$, where the theory degenerates to Euclidean harmonic analysis. An analogous characterisation is proved there by different means.

\section{Preliminaries}\label{sec:prelim}
This section recalls the notation, definitions, and prior results needed for the rest of the paper.

Let $\Hn=\C^n\times\R$ be the Heisenberg group with group law
\begin{align*}
(z,t)(w,s)=\Big(z+w,\ t+s+\tfrac12\operatorname{Im}(z\cdot\bar w)\Big).
\end{align*}
For $f\in L^1(\Hn), \ \mu\in\R^*$, write
\begin{align*}
f^\mu(z):=\int_\R f(z,t) e^{i\mu t} \,dt.
\end{align*}
 Let the Schr\"odinger representation $\pi_\lambda$ act on $L^2(\R^n)$, given explicitly by
\begin{align*}
\pi_\lambda(z,t)\varphi(\xi)=e^{i\lambda t}\,e^{i\lambda\left(x\cdot\xi+\frac12x\cdot y\right)}\varphi(\xi+y),\qquad \varphi\in L^2(\R^n),\ \xi\in\R^n,\ z=x+iy\in\C^n.
\end{align*}

The group Fourier transform of $f\in L^1(\Hn)$ is defined as 
\begin{align*}
\widehat{f}(\lambda)&:= \int_{\Hn} f(z,t) \ \pi_\lambda(z,t)\,dz dt.
\end{align*}

The Weyl transform of a function $g\in L^1(\C^n)$ is
\begin{align*}
W_\lambda(g) &:= \int_{\C^n}g(z)\,\pi_\lambda(z,0)\,dz.
\end{align*}

Write $f^*(z) := f(-z)$ for the involution on $\C^n$. The Weyl transform converts this involution into the operator adjoint.
\begin{align*}
W_\lambda(f^*) = W_\lambda(f)^*.
\end{align*}
For $\lambda\in\R^*$, the $\lambda$-twisted convolution of $f,g\in L^1(\C^n)$ is
\begin{align*}
(f*_\lambda g)(z)=\int_{\C^n}f(z-w)\,g(w)\,e^{i\frac\lambda2\operatorname{Im}(z\cdot\bar w)}\,dw.
\end{align*}
It is associative on $L^1(\C^n)$. The substitution $w\mapsto z-w$ gives the equivalent form
\begin{align*}
(f*_\lambda g)(z)=\int f(x)g(z-x)e^{-i\frac\lambda2\operatorname{Im}(z\cdot\bar x)}\,dx
\end{align*}
and the identity
\begin{align*}
f*_\mu g=g*_{-\mu}f.
\end{align*}
Since complex conjugation of $\C^n$ preserves radial functions and reverses the sign of $\operatorname{Im}(z\cdot\bar w)$, twisted convolution of radial functions does not depend on the sign of the parameter and is commutative.

The Weyl transform converts twisted convolution into composition of operators.
\begin{align*}
W_\lambda(f *_\lambda g) = W_\lambda(f)\,W_\lambda(g).
\end{align*}

For $f, g \in L^1(\Hn)$, the convolution $f * g$ is defined by
\begin{align*}
(f*g)(z,t) = \int_{\Hn} f\big((z,t)(-w,-s)\big)\,g(w,s)\,dw\,ds.
\end{align*}
The group Fourier transform converts convolution into composition of operators.
\begin{align*}
(f*g)\,\widehat{}\,(\lambda) = \widehat f(\lambda)\,\widehat g(\lambda).
\end{align*}

The involution is defined by $f^*(z,t) := f(-z,-t)$. The group Fourier transform converts the involution into the operator adjoint.
\begin{align*}
\widehat{f^*}(\lambda) = \widehat f(\lambda)^*.
\end{align*}

For $(z,t) \in \Hn$, right translation is defined by $R_{(z,t)}f(w,s) := f\big((w,s)(z,t)\big)$. The group Fourier transform satisfies
\begin{align*}
\big(R_{(z,t)}f\big)\,\widehat{}\,(\lambda) = \widehat f(\lambda)\,\pi_\lambda(z,t)^*.
\end{align*}

Let $L^\infty(\R^*, \mathcal B(L^2(\R^n)), d\mu)$ denote the space of essentially bounded $\mathcal B(L^2(\R^n))$-valued functions on $\R^*$, where $d\mu(\lambda) := (2\pi)^{-n-1}|\lambda|^n\,d\lambda$.

Convolution on $\Hn$ and twisted convolution on $\C^n$ are related by
\begin{align*}
(f*g)^\lambda = f^\lambda *_\lambda g^\lambda.
\end{align*}

For $\delta\ge-1/2$, let $L^\delta_k$ denote the Laguerre polynomials of type $\delta$, and for $\lambda\in\R^*$, $k\in\N$, define the Laguerre functions
\begin{align*}
\phik(z)=L^{n-1}_k\big(\tfrac12|\lambda||z|^2\big)\,e^{-\frac14|\lambda||z|^2}.
\end{align*}
These satisfy the uniform bound
\begin{align*}
|\phik(z)|\le\Lk, \qquad z\in\C^n
\end{align*}
(\cite[\S7.21]{Sz}).

The Heisenberg fan $\Omega\subset\R^2$ is
\begin{align*}
\Omega=\bigcup_{k=0}^\infty R_k\cup R_\infty,
\end{align*}
where
\begin{align*}
R_k=\{(\lambda,(2k+n)|\lambda|):\lambda\in\R^*\} \qquad\text{and}\qquad R_\infty=\{(0,\tau):\tau\ge0\}.
\end{align*} 

\begin{definition} (\cite{Th23})
\label{def:SFT}
For $f \in L^1 \cap L^2(\Hn)$, the Strichartz Fourier transform is defined as follows.

For $a = (\lambda,(2k+n)|\lambda|) \in R_k$,
\begin{align*}
\FS(f)(a,z) = \widehat f(a,z) := f^{-\lambda} *_{-\lambda} \varphi^{n-1}_{k,\lambda}(z).
\end{align*}

For $a = (0,\tau) \in R_\infty$ with $\tau > 0$,
\begin{align*}
\FS(f)(0,\tau,z) = \widehat f(0,\tau,z) := (n-1)!\,2^{n-1} \int_{\Hn} f(w,s)\,\frac{J_{n-1}\big(\sqrt\tau\,|z-w|\big)}{\big(\sqrt\tau\,|z-w|\big)^{n-1}}\,dw\,ds,
\end{align*}
where $J_{n-1}$ is the Bessel function of order $n-1$.

For $a = (0,0)$,
\begin{align*}
\FS(f)(0,0,z) = \widehat f(0,0,z) := \int_{\Hn} f(w,s)\,dw\,ds.
\end{align*}
\end{definition}

The special Hermite operator $L_{-\lambda}$ \cite{ThLect} has $L^2(\C^n)$-eigenspaces indexed by $j\in\N$, corresponding to eigenvalue $(2j+n)|\lambda|$. Let $P_j$ denote the orthogonal projection onto the $j$-th eigenspace. Then
\begin{align*}
u=\sum_{j=0}^\infty P_ju, \qquad u\in L^2(\C^n),
\end{align*}
and
\begin{align*}
u*_{-\lambda}\varphi^{n-1}_{j,\lambda}=(2\pi)^n|\lambda|^{-n}P_ju.
\end{align*}
The radial functions in the $j$-th eigenspace are the multiples of $\varphi^{n-1}_{j,\lambda}$, and radial tempered distributions on $\C^n$ have unique expansions $\sum_ja_j\varphi^{n-1}_{j,\lambda}$ converging in $\mathcal S'(\C^n)$ \cite{ThLect}. In particular, for $i,j\in\N$,
\begin{align}\label{eq:lag-orth}
\varphi^{n-1}_{i,\lambda}*_{-\lambda}\varphi^{n-1}_{j,\lambda}=\delta_{ij}(2\pi)^n|\lambda|^{-n}\varphi^{n-1}_{j,\lambda},
\end{align}
since $\varphi^{n-1}_{i,\lambda}$ lies in the $i$-th eigenspace. Hence its convolution on the right with $\varphi^{n-1}_{j,\lambda}$ is $(2\pi)^n|\lambda|^{-n}P_j\varphi^{n-1}_{i,\lambda}$. We note that radial twisted convolution is commutative.

The fact $\phik(0)=\Lk=\binom{k+n-1}{k}>0$ is used in Step~6 of
the proof of Theorem~\ref{thm:char} below, to fix the normalising
constant.

Finally, recall that $L_{-\lambda}$ is built from left-invariant vector fields on $\Hn$ via the Schr\"odinger representation, and that it commutes with twisted convolution on the left factor in the sense that, for Schwartz $\phi$ and suitable $\Theta$,
\begin{align}\label{eq:Lcomm}
L_{-\lambda}(\phi*_{-\lambda}\Theta)=(L_{-\lambda}\phi)*_{-\lambda}\Theta,
\end{align}
since $L_{-\lambda}$ is itself realised as convolution on the left by a (distributional) kernel built from $\phi$, and twisted convolution is associative \cite{ThLect}. This is used in Remark~\ref{rem:ellipticity} below.

The Strichartz Weyl transform was introduced in \cite{SWTpaper}, where its
basic structural theory was developed in detail. 

\begin{definition} (\cite{SWTpaper})
\label{def:SWT}
For $\lambda \in \R^*$ and $k \in \N$, the Strichartz Weyl transform is the operator
\begin{align} 
W_{S,\lambda,k}(g) := g *_{-\lambda} \varphi^{n-1}_{k,\lambda}, \qquad g \in L^1 \cap L^2(\C^n).\label{eq:WS}
\end{align}
When the index $k$ is understood, we write $W_{S,\lambda}(g)$ in place of $W_{S,\lambda,k}(g)$.
\end{definition}

Comparing Definition \ref{def:SFT} with Definition \ref{def:SWT} gives
\begin{align}
F_S(f)(a,z) = W_{S,\lambda}\big(f^{-\lambda}\big)(z), \qquad a = (\lambda,(2k+n)|\lambda|) \in R_k, \label{eq:FSWS}
\end{align}
for $f \in L^1 \cap L^2(\Hn)$, since $f^{-\lambda} \in L^1 \cap L^2(\C^n)$ for almost every $\lambda$ whenever $f \in L^1 \cap L^2(\Hn)$.

\begin{proposition} (\cite{SWTpaper})
\label{prop:SWT-proj}
For $\lambda \in \R^*$ and $k \in \N$, let $P_k$ denote the $k$-th spectral projection of $L_{-\lambda}$ on $L^2(\C^n)$. Then, for every $g \in L^1 \cap L^2(\C^n)$,
\begin{align*} 
W_{S,\lambda,k}(g) = (2\pi)^n|\lambda|^{-n} P_k g.
\end{align*}
\end{proposition}

The present paper needs
only four properties from that theory, stated and proved below as
Proposition~\ref{prop:props}, which serve as the hypotheses of the
characterisation theorem proved in Section~\ref{sec:charWS} (with (ii)
needed there only at the origin).

\begin{proposition}\label{prop:props}
The Strichartz Weyl transform satisfies the following properties.
\begin{enumerate}
\item[(i)] for every $u\in L^1(\C^n)$,
\begin{align*}
\WS(u*_{-\lambda}\phik)=(2\pi)^n|\lambda|^{-n}\,\WS(u),
\end{align*}
\item[(ii)] $\WS\big(\phik\big)=(2\pi)^n|\lambda|^{-n}\phik$. In particular,
\begin{align*}
\WS(\phik)(0)=(2\pi)^n|\lambda|^{-n}\phik(0),
\end{align*}
\item[(iii)] for every $c\in\C^n$ and $g\in L^1(\C^n)$,
\begin{align*}
\WS(\tau_cg)(w)=e^{i\frac\lambda2\operatorname{Im}(c\cdot\bar w)}\,\WS(M_cg)(w+c),
\end{align*}
\item[(iv)] for every $\sigma\in U(n)$ and $g\in L^1(\C^n)$,
\begin{align*}
\WS(\rho_\sigma g)(z)=\WS(g)(\sigma z),
\end{align*}
\end{enumerate}
where $\tau_cg(x):=e^{i\frac\lambda2\operatorname{Im}(x\cdot\bar c)}g(x+c)$, $M_cg(x):=e^{i\lambda\operatorname{Im}(x\cdot\bar c)}g(x)$, and $\rho_\sigma g(x):=g(\sigma x)$.
\end{proposition}
\begin{proof}
(i) By \eqref{eq:lag-orth} with $i=j=k$,
\begin{align}\label{eq:phik-self}
\phik*_{-\lambda}\phik=(2\pi)^n|\lambda|^{-n}\phik.
\end{align}
By associativity of twisted convolution and \eqref{eq:WS},
\begin{align*}
\WS(u*_{-\lambda}\phik)
&=u*_{-\lambda}\big(\phik*_{-\lambda}\phik\big)\\
&=u*_{-\lambda}\big((2\pi)^n|\lambda|^{-n}\phik\big)\\
&=(2\pi)^n|\lambda|^{-n}\,\WS(u).
\end{align*}

(ii) By Definition~\ref{def:SWT},
\begin{align*}
\WS(\phik) &= \phik *_{-\lambda} \phik \\
&= (2\pi)^n|\lambda|^{-n}\phik \qquad \text{(by \eqref{eq:phik-self})}.
\end{align*}

(iii) Note that,
\begin{align*}
\tau_cg(w-y)
&=e^{i\frac\lambda2\operatorname{Im}((w-y)\cdot\bar c)}g(w+c-y)\\
&=e^{i\frac\lambda2\operatorname{Im}(w\cdot\bar c)}e^{-i\frac\lambda2\operatorname{Im}(y\cdot\bar c)}g(w+c-y).
\end{align*}
Substituting $\varphi:=\phik$ in \eqref{eq:WS},
\begin{align}
\WS(\tau_cg)(w)&=\int_{\C^n}\tau_cg(w-y)\,\varphi(y)\,e^{-i\frac\lambda2\operatorname{Im}(w\cdot\bar y)}\,dy\notag\\
&=e^{i\frac\lambda2\operatorname{Im}(w\cdot\bar c)}\int_{\C^n} e^{-i\frac\lambda2\operatorname{Im}(y\cdot\bar c)}g(w+c-y)\,\varphi(y)\,e^{-i\frac\lambda2\operatorname{Im}(w\cdot\bar y)}\,dy.\label{eq:tauc}
\end{align}
On the other hand,
\begin{align*}
\WS(M_cg)(w+c)&=\int_{\C^n} M_cg(w+c-y)\,\varphi(y)\,e^{-i\frac\lambda2\operatorname{Im}((w+c)\cdot\bar y)}\,dy\\
&=\int_{\C^n} e^{i\lambda\operatorname{Im}(w\cdot\bar c)}e^{-i\lambda\operatorname{Im}(y\cdot\bar c)}g(w+c-y)\,\varphi(y)\,e^{-i\frac\lambda2\operatorname{Im}(w\cdot\bar y)}e^{i\frac\lambda2\operatorname{Im}(y\cdot\bar c)}\,dy\\
&=e^{i\lambda\operatorname{Im}(w\cdot\bar c)}\int_{\C^n} e^{-i\frac\lambda2\operatorname{Im}(y\cdot\bar c)}g(w+c-y)\,\varphi(y)\,e^{-i\frac\lambda2\operatorname{Im}(w\cdot\bar y)}\,dy\\
&=e^{-i\frac\lambda2\operatorname{Im}(w\cdot\bar c)}\cdot e^{i\frac\lambda2\operatorname{Im}(w\cdot\bar c)}\int_{\C^n} e^{-i\frac\lambda2\operatorname{Im}(y\cdot\bar c)}g(w+c-y)\,\varphi(y)\,e^{-i\frac\lambda2\operatorname{Im}(w\cdot\bar y)}\,dy\\
&=e^{-i\frac\lambda2\operatorname{Im}(c\cdot\bar w)}\,\WS(\tau_cg)(w).
\end{align*}
Combining the above with \eqref{eq:tauc}, we get
\begin{align*}
\WS(M_cg)(w+c)=e^{-i\frac\lambda2\operatorname{Im}(c\cdot\bar w)}\,\WS(\tau_cg)(w),
\end{align*}
which, since $e^{i\frac\lambda2\operatorname{Im}(c\cdot\bar w)}=e^{-i\frac\lambda2\operatorname{Im}(w\cdot\bar c)}$, is the stated identity (iii).

(iv) By \eqref{eq:WS},
\begin{align*}
\WS(\rho_\sigma g)(z)=\int_{\C^n} g(\sigma(z-y))\,\varphi(y)\,e^{-i\frac\lambda2\operatorname{Im}(z\cdot\bar y)}\,dy.
\end{align*}
Substituting $y'=\sigma y$, and using the fact that $\varphi$ is radial,
\begin{align*}
\WS(\rho_\sigma g)(z)
&=\int_{\C^n} g(\sigma(z-\sigma^{-1}y'))\,\varphi(\sigma^{-1}y')\,e^{-i\frac\lambda2\operatorname{Im}(z\cdot\overline{\sigma^{-1}y'})}\,dy'\\
&=\int_{\C^n} g(\sigma z-y')\,\varphi(y')\,e^{-i\frac\lambda2\operatorname{Im}((\sigma z)\cdot\bar y')}\,dy'\\
&=\WS(g)(\sigma z). \qedhere
\end{align*}
\end{proof}

\begin{theorem}[Lakshmi Lavanya-Thangavelu, \cite{LLT}]
\label{thm:LLT-weyl}
Let $T:\bigl(L^1(\C^n),*_\lambda\bigr)\to \mathcal B(L^2(\R^n))$ be a nonzero continuous $*$-homomorphism. Then there is a subspace $\mathcal H_\lambda$ of $L^2(\R^n)$ and a unitary representation $\rho_\lambda$ of $\Hn$ on $\mathcal H_\lambda$ such that
\begin{align*}
T(f) = \int_{\C^n} f(z)\,\rho_\lambda(z,0)\,dz, \quad \text{on } \mathcal H_\lambda,
\end{align*}
and there is a decomposition $L^2(\R^n) = \mathcal H_\lambda \oplus V^\lambda$, where
\begin{align*}
V^\lambda := \bigl\{v \in L^2(\R^n) : (Tf)(v) = 0 \text{ for all } f \in L^1(\C^n)\bigr\}.
\end{align*}
\end{theorem}

\begin{theorem}[Lakshmi Lavanya-Thangavelu, \cite{LLT}]
\label{thm:LLT-fourier}
Let $T:L^1(\Hn)\to L^\infty(\R^*, \mathcal B(L^2(\R^n)), d\mu)$ be a nonzero continuous linear map satisfying
\begin{enumerate}
  \item[(i)]
  \begin{align*}
  T(f^*)(\lambda) = T f(\lambda)^*, \quad \text{for all } \lambda \in \R^*,\ f \in L^1(\Hn),
  \end{align*}
  \item[(ii)]
  \begin{align*}
  T(f*g)(\lambda) = (Tf)(\lambda)\,(Tg)(\lambda), \quad \lambda \in \R^*,\ f,g \in L^1(\Hn),
  \end{align*}
  \item[(iii)]
  \begin{align*}
  T(R_{(z,t)} f)(\lambda) = (Tf)(\lambda)\,\pi_\lambda(z,t)^*, \quad \lambda \in \R^*,\ f \in L^1(\Hn),\ (z,t) \in \Hn.
  \end{align*}
\end{enumerate}
Then
\begin{align*}
(Tf)(\lambda) = \widehat f(\lambda), \quad \lambda \in A,\ f \in L^1(\Hn),
\end{align*}
where $A := \{\lambda \in \R^* : (Tf)(\lambda) \neq 0 \text{ for some } f \in L^1(\Hn)\}$.
\end{theorem}

Their proof proceeds in two stages. First, central translation covariance shows that $T(f)(\lambda)$ depends on $f$ only through $f^{\lambda}$. This reduces the problem to a map on $\C^n$. Second, the convolution property, now expressed for this reduced map, is shown to force it to be a Weyl transform $W_\lambda$. The special Hermite functions $\Phi_{\alpha\beta}$ multiply like matrix units, with $\Phi_{\alpha\beta}*_\lambda\Phi_{\mu\nu}$ being a multiple of $\delta_{\alpha\nu}\Phi_{\mu\beta}$. The images $Q_{\alpha\beta}$ of the normalised special Hermite functions inherit this relation, and an explicit representation can be reconstructed from it.

It is worth making the obstruction, and its resolution, explicit. The second
stage of the argument in~\cite{LLT} succeeds for a
structural reason. The special Hermite functions satisfy the matrix-unit
relation
\[
\Phi_{\alpha\beta} *_\lambda \Phi_{\mu\nu}
= \delta_{\alpha\nu}\,(2\pi)^{n/2}\,\Phi_{\mu\beta}.
\]
Any continuous $*$-homomorphism $T$ into $B(L^2(\mathbb{R}^n))$ carries this
relation forward to its images $Q_{\alpha\beta} := T(\Phi_{\alpha\beta})$.
This gives
\[
Q_{\alpha\beta} Q_{\mu\nu} = \delta_{\alpha\nu}\, Q_{\mu\beta}
\]
- a matrix-unit system inside $B(L^2(\mathbb{R}^n))$, from which an
explicit unitary representation can be reconstructed. This mechanism depends on the codomain being $B(L^2(\mathbb{R}^n))$ itself,
where composition is noncommutative and such a relation can hold.

For $\mathcal{F}_S$, the natural analogue $W_{S,\lambda}$ takes values in
$C(\mathbb{C}^n)$, a commutative, scalar-valued space. A relation of
the form $Q_{\alpha\beta} Q_{\mu\nu} = \delta_{\alpha\nu} Q_{\mu\beta}$
cannot be carried by scalar-valued images at all. Scalars commute, so no
system of nonzero scalars can satisfy a matrix-unit relation beyond the
trivial rank-one case. Correspondingly, $W_{S,\lambda}$ itself is not a
representation of the twisted convolution algebra but a (scaled) spectral
projection, idempotent on self-convolution, with no companion
off-diagonal structure to exploit. The two-stage strategy of \cite{LLT} therefore does not carry over here.
Its second stage has to be rebuilt from scratch.

The proof used here proceeds in two layers. Twisted translation
covariance together with $U(n)$-invariance already force the unknown kernel
to be radial (Steps~2-3 of \Cref{thm:char} below). Only then does a
single-function convolution-reproduction hypothesis suffice
to identify the radial dependence with the Laguerre function
$\varphi_{k,\lambda}^{n-1}$. This can be done either via the special Hermite expansion
(Step~5), or, as we show independently in \Cref{rem:ellipticity}, via the
ellipticity of the special Hermite operator $L_{-\lambda}$ alone. 

Throughout, $C(X)$ denotes continuous functions on $X$, with the topology of
uniform convergence on compact subsets. Under this topology, $C(X)$ is a
Fréchet space when $X$ is locally compact and $\sigma$-compact. Continuity
of linear maps into $C(X)$ is understood with respect to this topology.

We record here a measure-theoretic lemma used in Step~3 of the proofs of
\Cref{thm:char} and \Cref{thm:bdry} below. It upgrades an a.e.\ rotation-invariance
hypothesis, where the exceptional null set is allowed to depend on the
rotation, to a single a.e.\ radial representative.

\begin{lemma}\label{lem:radial}
Let $h \in L^\infty(\mathbb{C}^n)$. Suppose that for every $\sigma \in U(n)$
there is a Lebesgue-null set $N_\sigma \subset \mathbb{C}^n$ (depending on
$\sigma$) such that $h(\sigma x) = h(x)$ for every $x \notin N_\sigma$. Then
$h$ agrees with a radial function a.e.\ on $\mathbb{C}^n$.
\end{lemma}

\begin{proof}
\emph{Step 1 (exact invariance after mollifying).} Fix a radial
$\varphi \in C_c^\infty(\mathbb{C}^n)$, with $\varphi \ge 0$ and
$\int_{\mathbb{C}^n}\varphi = 1$. Fix
$\sigma \in U(n)$ and $y \in \mathbb{C}^n$. Then, using that $\sigma$ acts
on $\mathbb{C}^n \cong \mathbb{R}^{2n}$ as a Lebesgue-measure-preserving
orthogonal transformation and $\varphi(\sigma z) = \varphi(z)$ for every
$z$,
\begin{align}
(h*\varphi)(\sigma y) &= \int_{\mathbb{C}^n} h(x)\,\varphi(\sigma y - x)\,dx \notag\\
&= \int_{\mathbb{C}^n} h(\sigma x')\,\varphi(y-x')\,dx' \notag\\
&= (h_\sigma * \varphi)(y), \label{eq:step1a}
\end{align}
for every $y$, with no null sets involved. Here, $h_\sigma(x) := h(\sigma x)$.

Since $h_\sigma = h$ on the complement of the
null set $N_\sigma$,
\[
(h_\sigma * \varphi)(y) = (h*\varphi)(y), \qquad \text{for every } y \in \mathbb{C}^n.
\]

\emph{Step 2 (assembling over all of $U(n)$).} The conclusion of
Step~1 now holds for every $y$, for each
individually fixed $\sigma$, so we may let $\sigma$ range over all of
$U(n)$ with no obstruction from its uncountability. Thus $H := h*\varphi$
satisfies $H(\sigma y) = H(y)$ for every $\sigma \in U(n)$ and every $y$.
Since $U(n)$ acts transitively on each sphere in $\mathbb{C}^n \cong
\mathbb{R}^{2n}$, and $H$ is smooth (as $\varphi \in C_c^\infty$), $H$ is a radial function.

\emph{Step 3 (removing the mollifier).} Apply Steps~1-2 to a radial
approximate identity $\varphi_\varepsilon(x) :=
\varepsilon^{-2n}\varphi(x/\varepsilon) \to \delta_0$. Each
$h*\varphi_\varepsilon$ is radial, and $h*\varphi_\varepsilon \to h$
in $L^1_{\mathrm{loc}}(\mathbb{C}^n)$ as $\varepsilon \to 0$. Pass to a
subsequence $\varepsilon_m \to 0$ along which $h*\varphi_{\varepsilon_m} \to
h$ a.e.\ on $\mathbb{C}^n$. Write $h*\varphi_{\varepsilon_m}(x) =
r_m(|x|)$ for every $x$, since $h*\varphi_{\varepsilon_m}$ depends only on $|x|$.

In polar coordinates $\mathbb{C}^n \setminus \{0\} \cong (0,\infty) \times
S^{2n-1}$,
\[
dx = \rho^{2n-1}\,d\rho\,d\omega, \qquad x = \rho\,\omega, \quad \rho > 0,\ \omega \in S^{2n-1}.
\]
Since $h*\varphi_{\varepsilon_m} \to h$ a.e.\ on $\mathbb{C}^n$ (equivalently, a.e.\ with respect to $\rho^{2n-1}\,d\rho\,d\omega$), Fubini's theorem gives a set of full measure of $\rho > 0$ such that, for each such $\rho$,
\[
r_m(\rho) \longrightarrow h(\rho\,\omega) \quad \text{as } m \to \infty,
\quad \text{for a.e.\ } \omega \in S^{2n-1}.
\]

Fix such a $\rho$. The left side is a single sequence of numbers not
depending on $\omega$. If it failed to converge, the displayed statement
would fail for every $\omega$. This contradicts that it holds for a.e.\
$\omega$, as $S^{2n-1}$ has positive measure. Hence $r_m(\rho) \to
L(\rho)$ for some $L(\rho)$. By uniqueness of limits, every $\omega$ for
which the displayed convergence holds must then satisfy $h(\rho\,\omega) =
L(\rho)$. Thus $h(\rho\,\cdot\,) = L(\rho)$ for a.e.\ $\omega$, for a.e.\
$\rho$, that is, $h$ agrees a.e.\ on $\mathbb{C}^n$ with the radial
function $x \mapsto L(|x|)$.
\end{proof}

\section{Characterisation of the Strichartz Weyl Transform}\label{sec:charWS}

Here we characterise the Strichartz Weyl transform $\WS$ by four structural properties.

Recall the operators $\tau_c, M_c, \rho_\sigma$ introduced in Proposition~\ref{prop:props}.

\begin{theorem}\label{thm:char}
Fix $\lambda\in\R^*$ and $k\in\N$. Let $S:\big(L^1(\C^n),*_{-\lambda}\big)\to C(\C^n)$ be a continuous linear map satisfying
\begin{enumerate}
\item[(i)] for every $u\in L^1(\C^n)$,
\begin{align*}
S(u*_{-\lambda}\phik)=(2\pi)^n|\lambda|^{-n}\,S(u),
\end{align*}
\item[(ii)]
\begin{align*}
S\big(\phik\big)(0)=(2\pi)^n|\lambda|^{-n}\phik(0),
\end{align*}
\item[(iii)] for every $c\in\C^n$ and $g\in L^1(\C^n)$,
\begin{align*}
S(\tau_cg)(w)=e^{i\frac\lambda2\operatorname{Im}(c\cdot\bar w)}\,S(M_cg)(w+c),
\end{align*}
\item[(iv)] for every $\sigma\in U(n)$ and $g\in L^1(\C^n)$,
\begin{align*}
S(\rho_\sigma g)(z)=S(g)(\sigma z).
\end{align*}
\end{enumerate}
Then $S(u)=\WS(u)$ for every $u\in L^1(\C^n)$.
\end{theorem}

\begin{proof}
\emph{Step 1 (duality).} Since $S$ is a continuous linear map $L^1(\C^n)\to C(\C^n)$, for each fixed $w$ the functional $u\mapsto S(u)(w)$ is a bounded linear functional on $L^1(\C^n)$. By the duality $L^1(\C^n)^*=L^\infty(\C^n)$ there is a unique $h_w\in L^\infty(\C^n)$ with
\begin{align}
S(u)(w)=\int u(x)\,h_w(x)\,dx.\label{eq:Step1}
\end{align}

\emph{Step 2 ((iii) gives a twisted convolution kernel).} 
Writing out (iii) through Step~1,
\begin{align*}
\int_{\C^n} g(x-c)\,h_w(x-c)\,dx
&= e^{i\frac\lambda2\operatorname{Im}(c\cdot\bar w)}\int_{\C^n} g(x)\,e^{-i\frac\lambda2\operatorname{Im}((x-c)\cdot\bar c)}\,h_{w+c}(x-c)\,dx \\[4pt]
\Longrightarrow\quad \int_{\C^n} g(x)\Big[h_w(x-c) &- e^{i\frac\lambda2\operatorname{Im}(c\cdot\bar w)}e^{i\frac\lambda2\operatorname{Im}(x\cdot\bar c)}\,h_{w+c}(x-c)\Big]\,dx=0
\end{align*}
for every $g\in L^1(\C^n)$. Hence, for every $w,c\in\C^n$ and almost every $x$,
\begin{align*}
h_w(x-c)=e^{i\frac\lambda2\operatorname{Im}(c\cdot\bar w)}e^{i\frac\lambda2\operatorname{Im}(x\cdot\bar c)}\,h_{w+c}(x).
\end{align*}
(An $L^\infty$ function has no pointwise values, so we do not evaluate at a point. Instead we choose $c$ to reduce to $h_0$.) Taking $c=-w$ gives
\begin{align*}
h_w(x+w)=e^{-i\frac\lambda2\operatorname{Im}(x\cdot\bar w)}h_0(x)
\end{align*}
for a.e.\ $x$, that is,
\begin{align*}
h_w(y)=e^{-i\frac\lambda2\operatorname{Im}(y\cdot\bar w)}h_0(y-w).
\end{align*}
Put $\Theta(x):=h_0(-x)\in L^\infty(\C^n)$. Then
\begin{align*}
h_w(y)=e^{-i\frac\lambda2\operatorname{Im}(y\cdot\bar w)}\Theta(w-y).
\end{align*}
Applying \Cref{eq:Step1},
\begin{align}
S(u)(w) &= \int_{\C^n} u(y)\,e^{-i\frac\lambda2\operatorname{Im}(y\cdot\bar w)}\,\Theta(w-y)\,dy \notag\\
&= \big(u*_{-\lambda}\Theta\big)(w), \qquad u\in L^1(\C^n). \label{eq:Sconv}
\end{align}

\emph{Step 3 (rotation forces $\Theta$ radial).} By duality, (iv) gives, for each
fixed $\sigma \in U(n)$,
\[
h_z(\sigma^{-1}x) = h_{\sigma z}(x) \qquad \text{for a.e.\ } x \in
\mathbb{C}^n,
\]
with an exceptional null set that may depend on $\sigma$ (and on $z$). At
$z = 0$ this gives, for each fixed $\sigma \in U(n)$,
\[
h_0(\sigma^{-1}x) = h_0(x) \qquad \text{for a.e.\ } x \in \mathbb{C}^n.
\]
By Lemma~\ref{lem:radial} (applied with $\sigma$ in place of $\sigma^{-1}$,
which also ranges over all of $U(n)$), $h_0$, hence $\Theta$, agrees a.e.\
with a radial function.

\emph{Step 4 ((i) gives the eigen-identity for $\Theta$).} By \eqref{eq:Sconv} and associativity of twisted convolution,
\begin{align*}
S(u*_{-\lambda}\phik)=(u*_{-\lambda}\phik)*_{-\lambda}\Theta=u*_{-\lambda}\big(\phik*_{-\lambda}\Theta\big),
\end{align*}
so hypothesis (i) reads directly as
\begin{align}\label{eq:eig}
u*_{-\lambda}\big(\phik*_{-\lambda}\Theta\big)=(2\pi)^n|\lambda|^{-n}\,\big(u*_{-\lambda}\Theta\big)\quad\text{for every }u\in L^1(\C^n).
\end{align}

\emph{Step 5 (expansion of $\Theta$ and elimination of all but one coefficient).} Being bounded and radial (Step 3), $\Theta$ is a radial tempered distribution, with unique expansion $\Theta=\sum_ia_i\varphi^{n-1}_{i,\lambda}$ in $\mathcal S'(\C^n)$ (Section~\ref{sec:prelim}). Twisted convolution with the Schwartz function $\phik$ acts continuously on $\mathcal S'$, so by \eqref{eq:lag-orth}
\begin{align}
\phik*_{-\lambda}\Theta=\sum_ia_i\,\phik*_{-\lambda}\varphi^{n-1}_{i,\lambda}=(2\pi)^n|\lambda|^{-n}a_k\,\phik. \label{eq:thetak}
\end{align}
Apply \eqref{eq:eig} with $u=\varphi^{n-1}_{j,\lambda}$, $j\in\N$ (radial, Schwartz, hence integrable). Using \eqref{eq:thetak} and \eqref{eq:lag-orth} again,
\begin{align}
\varphi^{n-1}_{j,\lambda}*_{-\lambda}\big(\phik*_{-\lambda}\Theta\big)
&=(2\pi)^n|\lambda|^{-n}a_k\,\big(\varphi^{n-1}_{j,\lambda}*_{-\lambda}\phik\big)\notag\\
&=\delta_{jk}\,(2\pi)^{2n}|\lambda|^{-2n}a_k\,\phik. \label{eq:eig-left}
\end{align}
Expanding $\Theta$ in \eqref{eq:eig} and using \eqref{eq:lag-orth} once more,
\begin{align}
(2\pi)^n|\lambda|^{-n}\big(\varphi^{n-1}_{j,\lambda}*_{-\lambda}\Theta\big)
&=(2\pi)^n|\lambda|^{-n}\sum_ia_i\,\varphi^{n-1}_{j,\lambda}*_{-\lambda}\varphi^{n-1}_{i,\lambda}\notag\\
&=(2\pi)^{2n}|\lambda|^{-2n}a_j\,\varphi^{n-1}_{j,\lambda}. \label{eq:eig-right}
\end{align}
Combining \eqref{eq:eig}, \eqref{eq:eig-left} and \eqref{eq:eig-right} gives
\begin{align*}
\delta_{jk}a_k\phik=a_j\varphi^{n-1}_{j,\lambda}
\end{align*}
for every $j$. For $j\ne k$ this reads $a_j\varphi^{n-1}_{j,\lambda}=0$. Since $\varphi^{n-1}_{j,\lambda}(0)=L^{n-1}_j(0)\ne0$, $\varphi^{n-1}_{j,\lambda}\not\equiv0$, so $a_j=0$. Hence $\Theta=a_k\phik$.

\emph{Step 5 (expansion of $\Theta$ and elimination of all but one coefficient).} Being bounded and radial (Step 3), $\Theta$ is a radial tempered distribution, with unique expansion $\Theta=\sum_ia_i\varphi^{n-1}_{i,\lambda}$ in $\mathcal S'(\C^n)$ (Section~\ref{sec:prelim}). Twisted convolution with the Schwartz function $\phik$ acts continuously on $\mathcal S'$, so by \eqref{eq:lag-orth}
\begin{align*}
\phik*_{-\lambda}\Theta=\sum_ia_i\,\phik*_{-\lambda}\varphi^{n-1}_{i,\lambda}=(2\pi)^n|\lambda|^{-n}a_k\,\phik.
\end{align*}
Apply \eqref{eq:eig} with $u=\varphi^{n-1}_{j,\lambda}$, $j\in\N$ (radial, Schwartz, hence integrable). 

\emph{Step 6 (conclusion).} By \eqref{eq:Sconv} and \eqref{eq:lag-orth},
\begin{align*}
S\big(\phik\big)=a_k\,\phik*_{-\lambda}\phik=a_k(2\pi)^n|\lambda|^{-n}\phik.
\end{align*}
Evaluating at $0$ and comparing with (ii), using $\phik(0)=\Lk\ne0$, gives $a_k=1$. So $\Theta=\phik$, and \eqref{eq:Sconv} gives
\begin{align*}
S(u)=u*_{-\lambda}\phik=\WS(u).
\end{align*}
\end{proof}
Now we present an alternate proof to this result, which for most steps overlaps with the one given above.
\begin{remark}\label{rem:ellipticity}
Hypothesis (i) above asks only that the single function $\phik$ be reproduced correctly under twisted convolution. It does not quantify over an entire class of test functions. In this sense it is closer in spirit to the differentiation/shift framework of \cite{Khn,Ksft} than to a convolution-algebra hypothesis imposed on every $U(n)$-radial function. We now replace the special Hermite expansion argument of Step~5 above by a shorter argument from the ellipticity of the special Hermite operator $L_{-\lambda}$.

Steps 1-4 of the proof above apply unchanged. They give $S(u) = u *_{-\lambda} \Theta$ for a radial $\Theta \in L^\infty(\C^n)$, together with the identity \Cref{eq:eig}: for every $u \in L^1(\C^n)$,
\begin{align*}
u *_{-\lambda} \big(\phik *_{-\lambda} \Theta\big) = u *_{-\lambda} \big((2\pi)^n |\lambda|^{-n} \Theta\big).
\end{align*}

Fix $w \in \C^n$. By Step~1, $u \mapsto (u *_{-\lambda} \Phi)(w)$ is, for any fixed $\Phi \in L^\infty(\C^n)$, a bounded linear functional of $u \in L^1(\C^n)$. This functional is represented by a unique kernel in $L^\infty(\C^n)$. As a function of $x$, the kernel is built from $\Phi$ and the twisted-convolution phase, as in \Cref{eq:Sconv}. The identity above, namely \Cref{eq:eig}, says that the functionals obtained this way from $\Phi = \phik *_{-\lambda} \Theta$ and from $\Phi = (2\pi)^n |\lambda|^{-n} \Theta$ agree for every $u$. Hence the kernels representing them agree for a.e.\ $x$, for every $w$. Since the twisted-convolution phase is a nonvanishing factor common to both, this forces
\begin{align*}
\big(\phik *_{-\lambda} \Theta\big)(w - x) = (2\pi)^n |\lambda|^{-n} \Theta(w - x) \qquad \text{for a.e. } x, \text{ for every } w.
\end{align*}
Setting $x = 0$ and letting $w$ range over $\C^n$ gives the identity of functions on $\C^n$
\begin{align*}
\phik *_{-\lambda} \Theta = (2\pi)^n |\lambda|^{-n} \, \Theta.
\end{align*}

Apply $L_{-\lambda}$ to both sides. Recall that $L_{-\lambda}$ commutes with twisted convolution on the left factor, as noted in \Cref{eq:Lcomm} in Section~\ref{sec:prelim} (see \cite{ThLect}), and that $L_{-\lambda} \phik = (2k+n)|\lambda| \, \phik$. Together these give
\begin{align*}
L_{-\lambda} \Theta = (2k+n)|\lambda| \, \Theta.
\end{align*}
Since $L_{-\lambda}$ is elliptic with real-analytic coefficients, elliptic regularity for operators with real-analytic coefficients applies \cite{MN}. So $\Theta$, a priori only in $L^\infty(\C^n)$, is forced to be real-analytic.

Write $\Theta(z) = \theta(|z|)$, using that $\Theta$ is radial by Step~3. The eigenvalue equation becomes an ordinary differential equation in $r = |z|$. The substitution $x = \tfrac12 |\lambda| r^2$ turns it into
\begin{align*}
x y'' + n y' + \Bigl(\tfrac n2 + k - \tfrac x4\Bigr) y = 0, \qquad y(x) := \theta(r).
\end{align*}
Writing $y(x) = e^{-x/2} v(x)$ reduces this to the Laguerre equation of type $n-1$ and degree $k$ for $v$,
\begin{align*}
x v'' + (n - x) v' + k v = 0.
\end{align*}
This equation has a regular singular point at $x = 0$, with indicial roots $0$ and $-(n-1)$ \cite[Ch.~V]{Sz}. These are the same roots as for the equation in $y$. The coefficient of the term $-\tfrac{x}{4} y$ vanishes at $x = 0$, so this term does not enter the indicial equation. The root $0$ for $v$ gives the polynomial solution $v(x) = L_k^{n-1}(x)$, hence $y(x) = L_k^{n-1}(x) e^{-x/2}$, that is, $\theta(r) = \varphi_{k,\lambda}^{n-1}(z)$. This solution is real-analytic at $r = 0$. The second solution is not real-analytic at the origin. For $n \ge 2$ it has a power-type singularity at the origin, and for $n = 1$ a logarithmic one. Real-analyticity of $\Theta$ therefore excludes the second solution, so $\Theta$ is a scalar multiple of $\varphi_{k,\lambda}^{n-1}$.

Step~6 above now applies unchanged to fix this scalar at $1$.

This gives a second proof of Theorem~\ref{thm:char}, replacing the special Hermite expansion argument of Step~5 with an argument from ellipticity of $L_{-\lambda}$ alone. This is the same mechanism used for the single-ray uncertainty theorem of \cite{SWTpaper}.
\end{remark}

\section{Characterisation of the Strichartz Fourier Transform}\label{sec:charFS}

We now lift the characterisation of $\WS$ from Section~\ref{sec:charWS} to the Strichartz Fourier transform $\FS$ itself, on each ray of the Heisenberg fan.

For $(\zeta,\tau)\in\Hn$ and $\sigma\in U(n)$, let $R_{(\zeta,\tau)}$ denote right translation,
\begin{align*}
R_{(\zeta,\tau)}f(x,s):=f\big((x,s)(\zeta,\tau)\big).
\end{align*}
We write $M_\zeta$ and $\rho_\sigma$, the standard names for modulation and rotation, now for the corresponding operators on $\Hn$: $M_\zeta f(x,s):=e^{i\lambda\operatorname{Im}(x\cdot\bar\zeta)}f(x,s)$, $\rho_\sigma f(x,s):=f(\sigma x,s)$. In each case, the domain will be clear from context.

Consider the function
\begin{align*}
f_0(x,s):=\phik(x)\,e^{i\lambda s}e^{-s^2}\in L^1(\Hn),
\end{align*}
so that $f_0^{-\lambda}=\sqrt\pi\,\phik$.

First we summarise the required properties of the Strichartz Fourier transform.

\begin{proposition}\label{prop:FSprops}
For $a=(\lambda,(2k+n)|\lambda|)\in R_k$, the Strichartz Fourier transform $\FS$ satisfies
\begin{enumerate}
\item[(i)] for $f\in L^1(\Hn)$,
\begin{align*}
\FS(f*f_0)(a,z)=\sqrt\pi\,(2\pi)^n|\lambda|^{-n}\,\FS(f)(a,z),
\end{align*}
\item[(ii)]
\begin{align*}
\FS(f_0)(a,0)=\sqrt\pi\,(2\pi)^n|\lambda|^{-n}\phik(0),
\end{align*}
\item[(iii)] for $(\zeta,\tau)\in\Hn$,
\begin{align*}
\FS\big(R_{(\zeta,\tau)}f\big)(a,w)=e^{i\lambda\tau}e^{i\frac\lambda2\operatorname{Im}(\zeta\cdot\bar w)}\,\FS\big(M_\zeta f\big)(a,w+\zeta),
\end{align*}
\item[(iv)] for $\sigma\in U(n)$,
\begin{align*}
\FS(\rho_\sigma f)(a,w)=\FS(f)(a,\sigma w).
\end{align*}
\end{enumerate}
\end{proposition}

\begin{proof}
(i)  Recall
\begin{align*}
(f*g)^{-\lambda} = f^{-\lambda} *_{-\lambda} g^{-\lambda}.
\end{align*}
Apply this with $g=f_0$, so that $g^{-\lambda} = f_0^{-\lambda} = \sqrt\pi\,\phik$. By \Cref{eq:FSWS} and Proposition~\ref{prop:props}(i),
\begin{align*}
\FS(f*f_0)(a,z) &= \WS\big(f^{-\lambda}*_{-\lambda}\sqrt\pi\,\phik\big)(z)\\
&= \sqrt\pi\,(2\pi)^n|\lambda|^{-n}\WS\big(f^{-\lambda}\big)(z)\\
&= \sqrt\pi\,(2\pi)^n|\lambda|^{-n}\FS(f)(a,z),
\end{align*}
which is (i).

(ii) Since $f_0^{-\lambda}=\sqrt\pi\,\phik$, by \Cref{eq:FSWS} and Proposition \ref{prop:props}(ii),
\begin{align*}
\FS(f_0)(a,0)=\sqrt\pi\,\WS(\phik)(0)=\sqrt\pi\,(2\pi)^n|\lambda|^{-n}\phik(0).
\end{align*}

(iii) The group law gives $(x,s)(\zeta,\tau)=\big(x+\zeta,\,s+\tau+\frac12\operatorname{Im}(x\cdot\bar\zeta)\big)$. Substituting $s'=s+\tau+\frac12\operatorname{Im}(x\cdot\bar\zeta)$ gives
\begin{align*}
(R_{(\zeta,\tau)}f)^{-\lambda}(x)
&=e^{i\lambda\tau}e^{i\frac\lambda2\operatorname{Im}(x\cdot\bar\zeta)}f^{-\lambda}(x+\zeta)\\
&=e^{i\lambda\tau}\,\tau_\zeta\big(f^{-\lambda}\big)(x),
\end{align*}
Similarly, $(M_\zeta f)^{-\lambda}=M_\zeta(f^{-\lambda})$. By \Cref{eq:FSWS} and Proposition \ref{prop:props}(iii), this gives (iii).

(iv) Since $(\rho_\sigma f)^{-\lambda}=\rho_\sigma(f^{-\lambda})$, this follows from \Cref{eq:FSWS} and Proposition \ref{prop:props}(iv).
\end{proof}

This result characterises the Strichartz Fourier transform in terms of the properties discussed earlier.

\begin{theorem}\label{thm:FSchar}
Fix $\lambda\in\R^*$ and $k\in\N$, and put $a=(\lambda,(2k+n)|\lambda|)\in R_k$.\\
Let $T:L^1(\Hn)\to C(\Omega\times\C^n)$ be a continuous linear map such that $f\mapsto T(f)(a,\cdot)\in C(\C^n)$ satisfies
\begin{enumerate}
\item[(i)] for every $f\in L^1(\Hn)$,
\begin{align*}
T(f*f_0)(a,z)=\sqrt\pi\,(2\pi)^n|\lambda|^{-n}\,T(f)(a,z),
\end{align*}
\item[(ii)]
\begin{align*}
T(f_0)(a,0)=\sqrt\pi\,(2\pi)^n|\lambda|^{-n}\phik(0),
\end{align*}
\item[(iii)] for every $(\zeta,\tau)\in\Hn$,
\begin{align*}
T\big(R_{(\zeta,\tau)}f\big)(a,w)=e^{i\lambda\tau}e^{i\frac\lambda2\operatorname{Im}(\zeta\cdot\bar w)}\,T\big(M_\zeta f\big)(a,w+\zeta),
\end{align*}
\item[(iv)] for every $\sigma\in U(n)$,
\begin{align*}
T(\rho_\sigma f)(a,w)=T(f)(a,\sigma w).
\end{align*}
\end{enumerate}
Then $T(f)(a,z)=\FS(f)(a,z)$ for every $f\in L^1(\Hn)$ and every $z\in\C^n$.
\end{theorem}

\begin{proof}
\emph{Step A: Note that $T(f)(a,\cdot)$ depends on $f$ only through $f^{-\lambda}$.} Apply (iii) with $\zeta=0$. This gives
\begin{align*}
T(R_{(0,\tau)}f)(a,w)=e^{i\lambda\tau}T(f)(a,w),
\end{align*}
where $R_{(0,\tau)}f(x,s)=f(x,s+\tau)$. For $\eta\in L^1(\R)$ let
\begin{align*}
f_\eta:=\int_\R\eta(\tau)R_{(0,\tau)}f\,d\tau\in L^1(\Hn)
\end{align*}

By continuity of $T$,
\begin{align*}
T(f_\eta)(a,\cdot)=\eta^\lambda\,T(f)(a,\cdot)
\end{align*}
with $\eta^\lambda:=\int_\R\eta(\tau)e^{i\lambda\tau}\,d\tau$.

Take $\eta_R(\tau):=e^{-i\lambda\tau}\psi_R(\tau)$, with $\psi_R:=\frac1{2R}\mathbf 1_{[-R,R]}$, so that $\eta_R^\lambda=1$ and $T(f)(a,\cdot)=T(f_{\eta_R})(a,\cdot)$. If $f^{-\lambda}=0$, then
\begin{align*}
f_{\eta_R}(x,s)=e^{i\lambda s}(F_x*\psi_R)(s)
\end{align*}
with $F_x(u):=e^{-i\lambda u}f(x,u)$, and $\int_\R F_x=f^{-\lambda}(x)=0$.
Thus, $\|F*\psi_R\|_1\to0$ as $R\to\infty$. This is also true for $F=G'$ with $G\in C^1_c$, since then
\begin{align*}
\|G'*\psi_R\|_1\le\|G\|_1/R.
\end{align*}
The set of functions $G'$ with $G \in C^1_c(\mathbb{R})$ are $L^1$-dense in the zero-mean subspace $L^1_0(\mathbb{R})$. Combined with the uniform bound $\|F*\psi_R\|_1 \le \|F\|_1$ (valid for every $F\in L^1(\mathbb{R})$, since $\psi_R$ has mass $1$), this density gives the general case by a standard approximation argument.
By dominated convergence in $x$, $\|f_{\eta_R}\|_{L^1(\Hn)}\to0$. Hence $T(f)(a,\cdot)=T(f_{\eta_R})(a,\cdot)\to0$, i.e.\ $T(f)(a,\cdot)=0$.

Applying this to differences shows that $T(f)(a,\cdot)$ depends on $f$ only through $f^{-\lambda}$. So $S(g)(w):=T(f)(a,w)$ for any $f$ with $f^{-\lambda}=g$ is well defined on $L^1(\C^n)$. Every $g\in L^1(\C^n)$ arises as $f^{-\lambda}$ for some $f\in L^1(\Hn)$, namely $f:=g\otimes\eta$, where $(g\otimes\eta)(x,s):=g(x)\eta(s)$, for any $\eta\in L^1(\R)$ with $\int_\R e^{-i\lambda s}\eta(s)\,ds=1$ (such $\eta$ exists since the Fourier transform of an $L^1(\R)$ function need not vanish at $\lambda$). Moreover $g\mapsto g\otimes\eta$ is a bounded linear map $L^1(\C^n)\to L^1(\Hn)$ of norm $\|\eta\|_1$, so $S$ is linear and continuous.

\emph{Step B (reduction to Theorem \ref{thm:char}).} Specialise hypothesis (i) for $T$ to $g=f_0$, so that $g^{-\lambda}=f_0^{-\lambda}=\sqrt\pi\phik$. Using $(f*g)^{-\lambda}=f^{-\lambda}*_{-\lambda}g^{-\lambda}$, this reads
\begin{align*}
S\big(f^{-\lambda}*_{-\lambda}\sqrt\pi\phik\big)(z) &= T(f*f_0)(a,z)\\
&= \sqrt\pi\,(2\pi)^n|\lambda|^{-n}\,T(f)(a,z)\\
&= \sqrt\pi\,(2\pi)^n|\lambda|^{-n}\,S(f^{-\lambda})(z).
\end{align*}
Write $u=f^{-\lambda}$. By Step A, $u$ ranges over all of $L^1(\C^n)$ as $f$ ranges over $L^1(\Hn)$. This gives
\begin{align*}
S(u*_{-\lambda}\phik)=(2\pi)^n|\lambda|^{-n}\,S(u)\qquad\text{for every }u\in L^1(\C^n),
\end{align*}
which is hypothesis (i) of Theorem \ref{thm:char} for $S$.

Hypothesis (ii) for $T$ gives, via $f_0^{-\lambda}=\sqrt\pi\,\phik$ and linearity,
\begin{align*}
S(\phik)(0)=(2\pi)^n|\lambda|^{-n}\phik(0),
\end{align*}
which is hypothesis (ii) of Theorem \ref{thm:char}.

Recall
\begin{align*}
(R_{(\zeta,\tau)}f)^{-\lambda}(x)=e^{i\lambda\tau}\,\tau_\zeta\big(f^{-\lambda}\big)(x),\quad 
(M_\zeta f)^{-\lambda}=M_\zeta\big(f^{-\lambda}\big).
\end{align*}
Applying the above relations, hypothesis (iii) for $T$ with general $\zeta$ becomes hypothesis (iii) (of Theorem \ref{thm:char}) for $S$ (the factor $e^{i\lambda\tau}$ appears on both sides). By Theorem \ref{thm:char}, $S(g)=\WS(g)$ for every $g$. Hence, by \Cref{eq:FSWS},
\begin{align*}
T(f)(a,z)=S(f^{-\lambda})(z)=\WS(f^{-\lambda})(z)=\FS(f)(a,z).
\end{align*}
Hypothesis (iv) for $T$, via $(\rho_\sigma f)^{-\lambda}=\rho_\sigma(f^{-\lambda})$, becomes hypothesis (iv) for $S$. \qedhere
\end{proof}

\begin{remark}
The hypotheses (i)-(iv) are imposed, and the conclusion
established, only at the point $a \in R_k$, $k \in \mathbb{N}$, under
consideration. If they hold at every point of a ray $R_k$ (with $f_0$
depending on $a$), then $T = F_S$ on $R_k \times \mathbb{C}^n$. Nothing here
constrains $T$ on the boundary ray $R_\infty = \{(0,\tau) : \tau \ge 0\}$,
even though $T$'s stated codomain $C(\Omega \times \mathbb{C}^n)$ formally
includes it. The ray $R_\infty$ carries the degenerate
$\lambda = 0$ case, governed by the separate formula as in \Cref{def:SFT}, and none of (i)-(iv) as stated bears on it. We note that the
reduction to $f^{-\lambda}$ in Step~A of the proof above requires no
smoothness on $f$ and no spectral hypothesis on $T$, unlike the analogous
step in~\cite{Ksft}.
\end{remark}

\begin{remark}
By the reduction of Step~B, the alternative, ellipticity-based proof of Theorem \ref{thm:char} recorded in Remark \ref{rem:ellipticity} applies here unchanged, giving an alternate proof of Theorem \ref{thm:FSchar} as well.
\end{remark}

\section{The Strichartz Weyl Transform at the Boundary Ray}\label{sec:boundary}

At the boundary ray $R_\infty$, the twisted structure of Sections~\ref{sec:charWS}-\ref{sec:charFS} degenerates. We prove the analogous characterisations there by Euclidean methods.

\subsection{Arriving at the boundary}

Recall from Section~\ref{sec:prelim} that $R_\infty=\{(0,\tau):\tau\ge0\}$ is the boundary of the Heisenberg fan $\Omega$, where $\FS$ is denoted $\widehat f(0,\tau,z)$. We prove a characterisation for it mirroring Theorem \ref{thm:char}, with strictly weaker (untwisted) hypotheses. This is because the twisted structure underlying Section~\ref{sec:charWS} degenerates at $\lambda=0$.
We discuss the natural analogue $W_{S,\infty,\tau}$ of the Strichartz Weyl transform at $R_\infty$. This was obtained in \cite{SWTpaper}, and we include it here for completeness.

Fix $\tau>0$ and let $k\to\infty$, $\lambda\to0$ along $(2k+n)|\lambda|=\tau$. Two separate limits occur simultaneously.

\begin{enumerate}
\item The phase $e^{-i\lambda/2\operatorname{Im}(z\cdot\bar w)}$ defining $*_{-\lambda}$ tends to $1$. So twisted convolution degenerates to ordinary convolution $(u*v)(z)=\int_{\C^n}u(z-w)v(w)\,dw$ on $\C^n\cong\R^{2n}$.

\item The Mehler-Heine-type asymptotics for Laguerre polynomials give
\begin{align*}
\frac{\phik(z)}{\phik(0)}\longrightarrow K_\tau(z):=(n-1)!\,2^{n-1}\frac{J_{n-1}\big(\sqrt\tau|z|\big)}{\big(\sqrt\tau|z|\big)^{n-1}},\qquad k\to\infty,\ \lambda\to0,\ (2k+n)|\lambda|\to\tau,
\end{align*}
uniformly on compact sets \cite[Ch.~VIII]{Sz}. (See also \cite{Th23}, where they are used to establish continuity of $\FS$ at $R_\infty$.)
\end{enumerate}

Since $J_{n-1}(x)\sim(x/2)^{n-1}/(n-1)!$ as $x\to0$, we have $K_\tau(0)=1$ for every $\tau>0$. In fact
\begin{align*}
K_\tau(z)=\int_{S^{2n-1}}e^{i\sqrt\tau\,\omega\cdot z}\,d\omega
\end{align*}
(the Fourier transform of normalised surface measure on the sphere of radius $\sqrt\tau$ in $\R^{2n}$).

Combining these, as $k\to\infty$, $\lambda\to0$ along $(2k+n)|\lambda|=\tau$, the Strichartz Weyl transform $\WS(g)$ converges to $g*K_\tau$. This limit coincides with $W_{S,\infty,\tau}(g)$, the Strichartz Weyl transform at the boundary introduced in \cite{SWTpaper}, given by
\begin{align}\label{eq:WSinf}
W_{S,\infty,\tau}(g):=g*K_\tau.
\end{align}
Write $f^0(w):=\int_\R f(w,s)\,ds$, the $\mu=0$ case of $f^\mu$. Recall from Definition~\ref{def:SFT} that, for $a=(0,\tau,z)$ with $\tau>0$,
\begin{align*}
\FS(f)(0,\tau,z) := (n-1)!\,2^{n-1}\int_{\Hn} f(w,s)\,\frac{J_{n-1}(\sqrt\tau|z-w|)}{(\sqrt\tau|z-w|)^{n-1}}\,dw\,ds.
\end{align*}
The kernel does not depend on $s$, and equals $K_\tau(z-w)$. So
\begin{align*}
\FS(f)(0,\tau,z) &= \int_{\C^n}\Big(\int_\R f(w,s)\,ds\Big) K_\tau(z-w)\,dw\\
&= \int_{\C^n} f^0(w)\,K_\tau(z-w)\,dw\\
&= (f^0 * K_\tau)(z)\\
&= W_{S,\infty,\tau}(f^0)(z),
\end{align*}
which is the analogue of \Cref{eq:FSWS}.

Unlike $\WS$ ($\lambda\ne0$), which is (up to the scalar $(2\pi)^n|\lambda|^{-n}$) a projection onto a discrete eigenspace of $L_{-\lambda}$, $W_{S,\infty,\tau}$ is not idempotent. It coincides with the spectral density $dE_\tau/d\tau$ of the ordinary Laplacian $-\Delta$ on $\C^n\cong\R^{2n}$. The map $u\mapsto u*K_\tau$ restricts $\widehat u$ to the sphere $|\xi|=\sqrt\tau$, then reconstructs a function on $\C^n$ from this restriction. Equivalently, $u*K_\tau$ is the Hankel transform of $u$'s radial part, evaluated at $\sqrt\tau$, re-expanded in $z$. We write $H_\tau(u):=u*K_\tau$ for this operator.

\subsection{Characterisation of $W_{S,\infty,\tau}$}

\begin{theorem}[Characterisation of $W_{S,\infty,\tau}$]\label{thm:bdry}
Fix $\tau>0$. Let $S:\big(L^1(\C^n),*\big)\to C(\C^n)$ be a continuous linear map satisfying
\begin{enumerate}
\item[(i)] for every $u\in L^1(\C^n)$ and every $U(n)$-radial $v\in L^1(\C^n)$,
\begin{align*}
S(u*v)=S(v)(0)\,S(u),
\end{align*}
\item[(ii)] for two fixed, distinct positive reals $c_1,c_2$, writing $h_{c}(z):=e^{-c|z|^2}$,
\begin{align*}
S(h_{c_i})(0)=\Big(\frac\pi{c_i}\Big)^ne^{-\tau/(4c_i)},\qquad i=1,2,
\end{align*}
\item[(iii)] for every $c\in\C^n$ and $g\in L^1(\C^n)$,
\begin{align*}
S\big(g(\cdot+c)\big)(w)=S(g)(w+c),
\end{align*}
\item[(iv)] for every $\sigma\in U(n)$ and $g\in L^1(\C^n)$,
\begin{align*}
S(\rho_\sigma g)(z)=S(g)(\sigma z).
\end{align*}
\end{enumerate}
Then $S(u)=u*K_\tau=H_\tau(u)$ for every $u\in L^1(\C^n)$.
\end{theorem}

\begin{proof}
\emph{Step 1 (duality).} Since $S$ is continuous $L^1(\C^n)\to C(\C^n)$, for fixed $w$ the functional $u\mapsto S(u)(w)$ is a bounded linear functional on $L^1(\C^n)$. By $L^1(\C^n)^*\cong L^\infty(\C^n)$ there is a unique $h_w\in L^\infty(\C^n)$ with
\begin{align*}
S(u)(w)=\int_{\C^n} u(x)\,h_w(x)\,dx.
\end{align*}

\emph{Step 2 (translation covariance forces a convolution kernel).} Hypothesis (iii) gives
\begin{align*}
\int_{\C^n} g(x+c)\,h_w(x)\,dx = S\big(g(\cdot+c)\big)(w) = S(g)(w+c) = \int_{\C^n} g(x')\,h_{w+c}(x')\,dx'.
\end{align*}
Substituting $x'=x+c$ on the left,
\begin{align*}
\int_{\C^n} g(x')\,h_w(x'-c)\,dx' = \int_{\C^n} g(x')\,h_{w+c}(x')\,dx'.
\end{align*}
Since this holds for every $g\in L^1(\C^n)$, uniqueness of the $L^\infty$ kernel gives, for every $w,c$,
\begin{align*}
h_w(x'-c) = h_{w+c}(x') \qquad \text{for a.e. } x'.
\end{align*}
Take $c=-w$. (An $L^\infty$ function has no pointwise values, so we reduce to $h_0$ instead of evaluating at a point.) This gives
\begin{align*}
h_w(x'+w) = h_0(x'), \qquad \text{that is,} \qquad h_w(x) = h_0(x-w).
\end{align*}
Write $\Theta(x):=h_0(-x)$. Then $h_w(x)=\Theta(w-x)$, so $S(u)=u*\Theta$.

\emph{Step 3 (rotation forces $\Theta$ radial).} By duality, hypothesis (iv) gives,
for each fixed $\sigma \in U(n)$,
\[
h_z(\sigma^{-1}x) = h_{\sigma z}(x) \qquad \text{for a.e.\ } x \in
\mathbb{C}^n,
\]
with an exceptional null set that may depend on $\sigma$. At $z=0$ this
gives, for each fixed $\sigma \in U(n)$,
\[
h_0(\sigma^{-1}x) = h_0(x) \qquad \text{for a.e.\ } x \in \mathbb{C}^n.
\]
By Lemma~\ref{lem:radial}, $h_0$, hence $\Theta$, agrees a.e.\ with a
radial function.

\emph{Step 4 (radial convolution forces $\Theta=\gamma K_{\tau_0}$).} By associativity of convolution and hypothesis (i), for every radial $v\in L^1(\C^n)$ and every $u\in L^1(\C^n)$,
\begin{align*}
u*(v*\Theta)=S(u*v)=S(v)(0)\,u*\Theta.
\end{align*}
This implies
\begin{align}\label{eq:bdry-eig}
v*\Theta=S(v)(0)\,\Theta.
\end{align}
Hypothesis (ii) gives $S(h_{c_1})(0)=(\pi/c_1)^ne^{-\tau/(4c_1)}\ne0$. So $S\not\equiv0$, and by \Cref{eq:bdry-eig}, $\Theta\not\equiv0$.

Write $\widehat\Theta$ for the Fourier transform of $\Theta$ on $\widehat{\R^{2n}}\cong\C^n$. This is a tempered distribution, since $\Theta\in L^\infty$. Write $\widehat v$ for the ordinary Fourier transform of $v$. Then \Cref{eq:bdry-eig} becomes, for every radial $v\in L^1(\C^n)$,
\begin{align*}
\widehat v(\xi)\,\widehat\Theta(\xi)=S(v)(0)\,\widehat\Theta(\xi)\quad\text{in }\mathcal S'(\R^{2n}).
\end{align*}

Let $\chi(v):=S(v)(0)$. Associativity and \Cref{eq:bdry-eig} give
\begin{align*}
\chi(v_1*v_2)\Theta=\chi(v_1)\chi(v_2)\Theta
\end{align*}
for radial functions $v_1,v_2$. Since $\Theta\not\equiv0$, this gives $\chi(v_1*v_2)=\chi(v_1)\chi(v_2)$. Thus $\chi$ is a nonzero multiplicative linear functional on the commutative Banach algebra $\big(L^1_{rad}(\C^n),*\big)$. By the Gelfand theory, its characters are the maps $v\mapsto\widehat v(\rho)$ for $\rho\ge0$, that is, point evaluations of the radial Euclidean Fourier transform on $\R^{2n}$ \cite{SW}. Hence
\begin{align*}
\chi(v)=\widehat v(\rho_0)
\end{align*}
for every radial function $v$ and some fixed $\rho_0\ge0$.

Put $m(\xi):=|\xi|^2-\rho_0^2$. Fix the radial Schwartz function $v_0$ with $\widehat v_0(\xi):=m(\xi)e^{-m(\xi)^2}$. It is smooth, rapidly decreasing, and vanishing on the sphere $|\xi|=\rho_0$. In particular $\chi(v_0)=\widehat v_0(\rho_0)=0$. Then \Cref{eq:bdry-eig} at $v=v_0$ gives
\begin{align*}
\widehat v_0(\xi)\,\widehat\Theta(\xi)=0 \quad \text{in } \mathcal S'(\R^{2n}).
\end{align*}
The function $\widehat v_0$ is smooth and nonvanishing off the sphere $\{|\xi|=\rho_0\}$. So, for any test function $\phi$ supported away from this sphere, $\phi/\widehat v_0$ remains in $C_c^\infty(\R^{2n})$. Since $\widehat v_0\,\widehat\Theta=0$ as a distribution, pairing both sides against $\phi/\widehat v_0$ gives $\widehat\Theta(\phi)=0$, because $\widehat v_0\cdot(\phi/\widehat v_0)=\phi$ on the support of $\phi$. Hence
\begin{align*}
\operatorname{supp}\widehat\Theta\subset\{\xi:|\xi|=\rho_0\}.
\end{align*}

This support is compact and hence allows us to exclude derivatives of the surface measure. Choose $\kappa\in C_c^\infty(\R^{2n})$ with $\kappa=1$ near the sphere. The function $\kappa e^{m^2}$ is smooth, and $\kappa e^{m^2}\widehat v_0=\kappa m$.

Multiplying $\widehat v_0\widehat\Theta=0$ by $\kappa e^{m^2}$ gives
\begin{align*}
\big(\kappa e^{m^2}\widehat v_0\big)\widehat\Theta=0, \qquad\text{that is,}\qquad \kappa m\,\widehat\Theta=0.
\end{align*}
Since $\kappa\widehat\Theta=\widehat\Theta$ on the support of $\widehat\Theta$, this gives $m\widehat\Theta=0$, equivalently, $(\Delta+\rho_0^2)\Theta=0$.

By ellipticity, $\Theta$ is smooth. Writing $\Theta(z)=\theta(|z|)$, this becomes the radial form of Bessel's equation of order $n-1$,
\begin{align*}
\theta''+\frac{2n-1}r\theta'+\rho_0^2\theta=0.
\end{align*}
For $\rho_0>0$, the equation has two linearly independent solutions, $J_{n-1}(\rho_0 r)/r^{n-1}$ and $Y_{n-1}(\rho_0 r)/r^{n-1}$. At $r=0$, the first solution is regular, and the second is singular (\cite{Wat}). For $\rho_0=0$, the two solutions are $1$ and $r^{2-2n}$, or $\log r$ if $n=1$. Since $\Theta\in L^\infty(\C^n)$ is bounded near the origin, the singular solution is excluded. So $\theta$ is a scalar multiple of the regular solution,
\begin{align*}
\Theta=\gamma\,K_{\tau_0},\qquad\tau_0:=\rho_0^2\ge0.
\end{align*}

\emph{Step 5 (normalisation pins down $\gamma$ and $\tau_0$).} We compute $(h_c*K_{\tau_0})(0)$.
\begin{align*}
(h_c*K_{\tau_0})(0) &= \int_{\C^n} e^{-c|w|^2}\,K_{\tau_0}(w)\,dw\\
&= \int_{S^{2n-1}}\left(\int_{\C^n} e^{-c|w|^2}\,e^{i\sqrt{\tau_0}\,\omega\cdot w}\,dw\right)d\omega\\
&= \int_{S^{2n-1}} (\pi/c)^n e^{-\tau_0/(4c)}\,d\omega\\
&= (\pi/c)^n e^{-\tau_0/(4c)}
\end{align*}
for every $\tau_0\ge0$, $c>0$. Hypothesis (ii), applied at $c=c_1$ and $c=c_2$, gives
\begin{align*}
\gamma\,e^{-\tau_0/(4c_1)}=e^{-\tau/(4c_1)},\qquad \gamma\,e^{-\tau_0/(4c_2)}=e^{-\tau/(4c_2)}.
\end{align*}
Since $\tau_0,\tau,c_1,c_2$ are real, each equation forces $\gamma=e^{(\tau_0-\tau)/(4c_1)}=e^{(\tau_0-\tau)/(4c_2)}$. By injectivity of the exponential function on $\mathbb{R}$,
\begin{align*}
\frac{\tau_0-\tau}{4c_1}=\frac{\tau_0-\tau}{4c_2},
\end{align*}
and since $c_1\ne c_2$ this forces $\tau_0=\tau$, whence $\gamma=1$. Hence $\Theta=K_\tau$, and by Step 2,
\begin{align*}
S(u)=u*K_\tau=H_\tau(u)
\end{align*}
for every $u\in L^1(\C^n)$.
\end{proof}

\begin{remark}\label{rem:bdry}
Theorem \ref{thm:bdry} is not a limiting case of Theorem \ref{thm:char}. At $\lambda=0$ the projective phase in hypotheses (iii)/(iv) of Theorem \ref{thm:char} trivialises, so (iii) and (iv) above are (untwisted) translation and rotation covariance. The proof correspondingly replaces the Laguerre-expansion argument of Section~\ref{sec:charWS} with facts from Euclidean harmonic analysis on $\C^n\cong\R^{2n}$:
\begin{enumerate}
\item Translation invariance forces a convolution kernel.
\item The Gelfand theory of the radial convolution algebra identifies its possible eigenvalues.
\item Elementary distribution theory, rather than any Tauberian theorem, localises $\widehat\Theta$ to a single sphere.
\item The regularity of solutions of the radial Bessel equation, together with the Gaussian computation of Step 5, completely determines the kernel.
\end{enumerate}
Hypothesis (ii) also differs from its counterpart in Theorem~\ref{thm:char}. The boundary case determines a point of a continuous spectrum rather than a discrete eigenvalue. So a single normalising value constrains only one relation between the two unknowns $\gamma$ and $\tau_0$. Two independent values are needed to determine both, in contrast with Theorem~\ref{thm:char}'s single normalisation, which fixes only one unknown scalar. Unlike Theorem~\ref{thm:char}, hypothesis (i) here is not weakened to a single test function. The Gelfand-theoretic argument of Step~4 requires testing against every radial $v$, since it is the full multiplicativity of $\chi(\cdot)$ on $L^1_{\mathrm{rad}}(\C^n)$ that identifies $\rho_0$.
\end{remark}

\subsection{Characterisation of $\FS$ at $R_\infty$}

For $(\zeta,\theta)\in\Hn$ and $\sigma\in U(n)$, write 
\begin{align*}
R_{(\zeta,\theta)}f(x,s):=f\big((x,s)(\zeta,\theta)\big), \qquad \rho_\sigma f(x,s):=f(\sigma x,s).
\end{align*}
Fix $c_0>0$ and set
\begin{align*}
f_{c_0}(x,s):=h_{c_0}(x)\eta(s)\in L^1(\Hn),
\end{align*}
for any $\eta\in L^1(\R)$ with $\int_\R\eta=1$, so that $f_{c_0}^0=h_{c_0}$.

\begin{proposition}\label{prop:bdryFS}
With notations as defined above, the Strichartz Fourier transform $\FS$ satisfies the following, at $(0,\tau)$
\begin{enumerate}
\item[(i)] for $f\in L^1(\Hn)$, $g\in L^1(\Hn/U(n))$ radial,
\begin{align*}
\FS(f*g)(0,\tau,z)=\FS(g)(0,\tau,0)\,\FS(f)(0,\tau,z),
\end{align*}
\item[(ii)] for every $c_0>0$,
\begin{align*}
\FS(f_{c_0})(0,\tau,0)=\big(\pi/c_0\big)^ne^{-\tau/(4c_0)},
\end{align*}
\item[(iii)] for every $(\zeta,\theta)\in\Hn$,
\begin{align*}
\FS\big(R_{(\zeta,\theta)}f\big)(0,\tau,w)=\FS(f)(0,\tau,w+\zeta),
\end{align*}
\item[(iv)] for every $\sigma\in U(n)$,
\begin{align*}
\FS(\rho_\sigma f)(0,\tau,w)=\FS(f)(0,\tau,\sigma w).
\end{align*}
\end{enumerate}
\end{proposition}

\begin{proof}
Write $\FS(f)(0,\tau,z)=(f^0*K_\tau)(z)$, as established in Section~\ref{sec:boundary}.

(i) A direct computation gives $(f*g)^0=f^0*g^0$. Since $g$ is $U(n)$-radial, so is $g^0$. For any radial $v\in L^1(\C^n)$, write $\widehat v$ for the ordinary Fourier transform on $\R^{2n}\cong\C^n$. Then $\widehat v$ is also radial and hence it is constant on the sphere $|\xi|=\sqrt\tau$, where $\widehat K_\tau$ is supported. So
\begin{align*}
v*K_\tau=\widehat v(\sqrt\tau)K_\tau.
\end{align*}
Evaluating at $0$ and using $K_\tau(0)=1$ gives $\widehat v(\sqrt\tau)=(v*K_\tau)(0)$. That is,
\begin{align*}
v*K_\tau=(v*K_\tau)(0)\,K_\tau
\end{align*}
for radial $v$. Applying this with $v=g^0$, and using associativity of ordinary convolution gives
\begin{align*}
\FS(f*g)(0,\tau,z)
&=f^0*(g^0*K_\tau)(z)\\
&=(g^0*K_\tau)(0)\,(f^0*K_\tau)(z)\\
&=\FS(g)(0,\tau,0)\,\FS(f)(0,\tau,z).
\end{align*}

(ii) Since $f_{c_0}^0=h_{c_0}$, the computation of Step~5 above gives
\begin{align*}
\FS(f_{c_0})(0,\tau,0)=(h_{c_0}*K_\tau)(0)=(\pi/c_0)^ne^{-\tau/(4c_0)}.
\end{align*}

(iii) Substituting $s'=s+\theta+\tfrac12\operatorname{Im}(x\cdot\bar\zeta)$,
\begin{align*}
(R_{(\zeta,\theta)}f)^0(x) &= \int_\R f\big(x+\zeta,\, s+\theta+\tfrac12\operatorname{Im}(x\cdot\bar\zeta)\big)\,ds\\
&= \int_\R f(x+\zeta,s')\,ds'\\
&= f^0(x+\zeta).
\end{align*}
Hence $\FS\big(R_{(\zeta,\theta)}f\big)(0,\tau,w) = \big(f^0(\cdot+\zeta)*K_\tau\big)(w)$. Substituting $x'=x+\zeta$ in this convolution integral gives
\begin{align*}
(f^0*K_\tau)(w+\zeta) = \FS(f)(0,\tau,w+\zeta).
\end{align*}

(iv) Since $(\rho_\sigma f)^0=\rho_\sigma(f^0)$ and $K_\tau$ is radial, substituting $x=\sigma x'$ in the convolution integral defining $\FS(\rho_\sigma f)(0,\tau,w)$ gives
\begin{align*}
\FS(f)(0,\tau,\sigma w)=\FS(\rho_\sigma f)(0,\tau,w).
\end{align*}
\end{proof}

\begin{theorem}[Characterisation of $\FS$ at $R_\infty$]\label{thm:FSbdry}
Fix $\tau>0$. Let $T:L^1(\Hn)\to C(\C^n)$ be a continuous linear map satisfying, for every $f\in L^1(\Hn)$
\begin{enumerate}
\item[(i)] for every $U(n)$-radial $g\in L^1(\Hn/U(n))$,
\begin{align*}
T(f*g)(z)=T(g)(0)\,T(f)(z),
\end{align*}
\item[(ii)] for two fixed, distinct positive reals $c_1,c_2$,
\begin{align*}
T(f_{c_i})(0)=\big(\pi/c_i\big)^ne^{-\tau/(4c_i)},\qquad i=1,2,
\end{align*}
\item[(iii)] for every $(\zeta,\theta)\in\Hn$,
\begin{align*}
T\big(R_{(\zeta,\theta)}f\big)(w)=T(f)(w+\zeta),
\end{align*}
\item[(iv)] for every $\sigma\in U(n)$,
\begin{align*}
T(\rho_\sigma f)(w)=T(f)(\sigma w).
\end{align*}
\end{enumerate}
Then $T(f)(z)=\FS(f)(0,\tau,z)$ for every $f\in L^1(\Hn)$.
\end{theorem}

\begin{proof}
Applying (iii) with $\zeta=0$, we get
\begin{align*}
T(R_{(0,\theta)}f)(w)=T(f)(w)
\end{align*}
for every $\theta\in\R$. Unlike the case $a\in R_k$, no phase factor appears in this identity. Fix $\eta\in L^1(\R)$ with $\int_\R\eta=1$ and set
\begin{align*}
f_\eta(x,s):=\int_\R\eta(\theta)f(x,s+\theta)\,d\theta.
\end{align*}
Since $T$ is continuous and linear, this gives $T(f_\eta)(w)=T(f)(w)$ for every such $\eta$. A direct computation gives $f_\eta^0=f^0$ for every such $\eta$.

The argument of Step~A in the proof of Theorem~\ref{thm:FSchar} applies without change, with $\eta_R:=\frac1{2R}\mathbf 1_{[-R,R]}$ in place of $\eta_R$ used there and $f^0$ in place of $f^{-\lambda}$. If $h\in L^1(\Hn)$ has $h^0=0$, then $T(h)(w)=0$. Applying this to $h=f_1-f_2$ whenever $f_1^0=f_2^0$ shows $T(f)(w)$ depends on $f$ only through $f^0$.

So $S(g)(w):=T(f)(w)$ for any $f$ with $f^0=g$ is a well-defined continuous linear map on $L^1(\C^n)$. Every $g\in L^1(\C^n)$ arises as $f^0$ for some $f\in L^1(\Hn)$, e.g.\ $f(x,s):=g(x)\eta(s)$ for any $\eta\in L^1(\R)$ with $\int_\R\eta=1$. Continuity of $S$ follows from that of $T$ as in the proof of Theorem~\ref{thm:FSchar}.

Hypothesis (i) for $T$ gives, via $(f*g)^0=f^0*g^0$, hypothesis (i) of Theorem~\ref{thm:bdry} for $S$.

Hypothesis (ii) for $T$ translates, via $f_{c_i}^0=h_{c_i}$ ($i=1,2$), into hypothesis (ii) for $S$.

Hypothesis (iii) for $T$ with general $\zeta$ reduces, via
\begin{align*}
(R_{(\zeta,\theta)}f)^0(x)=f^0(x+\zeta),
\end{align*}
to hypothesis (iii) for $S$.

By Theorem~\ref{thm:bdry}, $S(g)=g*K_\tau$ for every $g$. Hence
\begin{align*}
T(f)(z)=S(f^0)(z)=(f^0*K_\tau)(z)=\FS(f)(0,\tau,z).
\end{align*}

Hypothesis (iv) for $T$ corresponds, via $(\rho_\sigma f)^0=\rho_\sigma(f^0)$, to hypothesis (iv) for $S$.
\end{proof}

\section{Conclusion}

We characterised $\WS$, the Strichartz Weyl transform, by four properties. These are a convolution-reproduction property at the single test function $\phik$, a twisted translation-covariance property, a $U(n)$-invariance property, and a single normalisation (Theorem~\ref{thm:char}). Lifting this via full translation covariance on $\Hn$ gives a characterisation of $\FS$ itself on each ray $R_k$, $k\in\N$, of the Heisenberg fan, again by a convolution-reproduction property at a single test function, $f_0$ (Theorem~\ref{thm:FSchar}).

As recorded in Remark~\ref{rem:ellipticity}, this single-function hypothesis also admits a proof via the ellipticity of $L_{-\lambda}$, in place of the special Hermite expansion argument of Section~\ref{sec:charWS}. By the reduction underlying Theorem~\ref{thm:FSchar}, the same alternative proof carries over unchanged to $\FS$.

At the boundary ray $R_\infty$, the twisted structure degenerates. In this case, $\WS$ and $\FS$ reduce to the spectral density of the ordinary Laplacian on $\C^n\cong\R^{2n}$ and to the Hankel transform, respectively. We proved an analogous characterisation there by Euclidean methods (Section~\ref{sec:boundary}), so that $\FS$ is characterised at every point of $\Omega$. At the boundary ray, the radial convolution hypothesis requires the full family of radial test functions, rather than a single one (Remark~\ref{rem:bdry}).

The characterisation of $\WS$ required a different route from the classical Weyl-transform-based argument of \cite{LLT}. Since $\WS$ is a projection rather than a faithful algebra representation, the matrix-unit technique driving that argument does not apply. Twisted translation-covariance and $U(n)$-invariance, together with the single-function convolution-reproduction property, close this gap and determine $\WS$ completely. The same pattern recurs at the boundary ray. There, translation covariance again reduces the problem to a convolution-algebra condition, closed by a radial eigenvalue argument, carried out by Euclidean methods instead.


\begin{thebibliography}{99}

\bibitem{ABM1} S. Achar, M.H.A. Biswas, R. Manna, \textit{Characterization of Fourier transform of Hilbert space valued functions,} arXiv:2610.36957.

\bibitem{ABM2} S. Achar, M.H.A. Biswas, R. Manna, \textit{Characterization of the group Fourier transform of Hilbert space valued functions,} arXiv:2611.03267.

\bibitem{AM} S. Achar, R. Manna, \textit{Characterization of Integral Transforms for Hilbertspace-Valued Functions,} arXiv:2611.04698.

\bibitem{AAM}
S.~Alesker, S.~Artstein-Avidan, and V.~Milman, \emph{A characterization of the Fourier transform and related topics}, Amer.\ Math.\ Soc.\ Transl.\ Ser.~2, Vol.~226, Amer.\ Math.\ Soc., Providence, RI, 2009, pp.~11-26.

\bibitem{DDG}
A.~Dabra, A.~Dasgupta, and P.~Gulia, \emph{Uncertainty principles for the Strichartz Fourier transform on the Heisenberg group}, Bull.\ Lond.\ Math.\ Soc.\ \textbf{58} (2026), e70424, doi:10.1112/blms.70424.

\bibitem{DGPT}
A.~Dasgupta, P.~Gulia, S.~Pusti, and S.~Thangavelu, \emph{Weighted norm inequalities for the Strichartz Fourier transform on the Heisenberg group}, arXiv:2603.02144, 2026.

\bibitem{Jaming}
P.~Jaming, \emph{A characterization of Fourier transforms}, Colloq.\ Math.\ \textbf{118} (2010), no.~2, 569-580.

\bibitem{K1}
M.~Krukowski, \emph{Integral transforms characterized by convolution}, Results Math. \textbf{78} (2023), no.~4, Paper No.~120, 12~pp.

\bibitem{K2}
M.~Krukowski, \emph{Characterizing the Fourier transform by its properties}, Pacific J. Math. \textbf{329} (2024), no.~2, 217-232.

\bibitem{Khn}
R.~Lakshmi Lavanya, \emph{Differentiation and Shift Property Characterisations of the Group Fourier Transform on the Heisenberg Group}, arXiv:2610.40508.

\bibitem{Ksft}
R.~Lakshmi Lavanya, \emph{Some Characterisations of the Strichartz Fourier Transform on the Heisenberg Group}, arXiv:2611.01128.

\bibitem{SWTpaper}
R.~Lakshmi Lavanya, \emph{On a Scalar-valued Weyl transform associated to the Heisenberg group}, arXiv:2611.36402.

\bibitem{LLT}
R.~Lakshmi Lavanya and S.~Thangavelu, \emph{A characterisation of the Fourier transform on the Heisenberg group}, Ann.\ Funct.\ Anal.\ \textbf{3} (2012), no.~1, 109-120, doi:10.15352/afa/1399900028.

\bibitem{MN}
C.~B.~Morrey, Jr.\ and L.~Nirenberg, \emph{On the analyticity of the solutions of linear elliptic systems of partial differential equations}, Comm.\ Pure Appl.\ Math.\ \textbf{10} (1957), 271-290.

\bibitem{SW}
E.~M.~Stein and G.~Weiss, \emph{Introduction to Fourier Analysis on Euclidean Spaces}, Princeton University Press, Princeton, NJ, 1971.


\bibitem{KS}
N.~Shravan Kumar and S.~Sivananthan, \emph{Characterisation of the Fourier transform on compact groups}, Bull. Aust. Math. Soc. \textbf{93} (2016), no.~3, 467-472.

\bibitem{Sz}
G.~Szeg\H{o}, \emph{Orthogonal Polynomials}, 4th ed., American Mathematical Society Colloquium Publications, Vol.~23, American Mathematical Society, Providence, RI, 1975.

\bibitem{ThLect}
S.~Thangavelu, \emph{Lectures on Hermite and Laguerre Expansions}, Mathematical Notes 42, Princeton University Press, Princeton, NJ, 1993.

\bibitem{ThHA}
S.~Thangavelu, \emph{Harmonic Analysis on the Heisenberg Group}, Progress in Mathematics 159, Birkh\"auser, Boston, MA, 1998.


\bibitem{Th23}
S.~Thangavelu, \emph{A scalar-valued Fourier transform for the Heisenberg group}, in: From Classical Analysis to Analysis on Fractals: A Tribute to Robert Strichartz, Vol.~1, Applied and Numerical Harmonic Analysis, Birkh\"auser, Cham, 2023, pp.~137-163, https://doi.org/10.1007/978-3-031-37800-3\_7.

\bibitem{Wat}
G.~N.~Watson, \emph{A Treatise on the Theory of Bessel Functions}, 2nd ed., Cambridge University Press, Cambridge, 1944.

\end{thebibliography}
\end{document}